\documentclass{amsart}

\usepackage{amssymb}

\usepackage{epsfig}         

\newtheorem{thm}{Theorem}[section]
\newtheorem{prop}[thm]{Proposition}

\newtheorem{cor}[thm]{Corollary}

\theoremstyle{definition}

\theoremstyle{remark}

\newtheorem{remark}[thm]{Remark}

\numberwithin{equation}{section}
\newcommand{\N}{\mathbb{N}}  
\newcommand{\R}{\mathbb{R}}  
\newcommand{\C}{\mathbb{C}}  

\begin{document}

\title{\sc Delay differential equations with differentiable solution operators on open domains in $C((-\infty,0],\R^n)$,
and processes for Volterra integro-differential equations}

\author{Hans-Otto Walther}

\address{Mathematisches Institut, Universit\"{a}t Gie{\ss}en,
Arndtstr. 2, D 35392 Gie{\ss}en, Germany. E-mail {\tt
Hans-Otto.Walther@math.uni-giessen.de}, phone ++49-151-58546608,
fax ++49-641-9932029}

\date{January 26, 2018}

\begin{abstract}

For autonomous delay differential equations $x'(t)=f(x_t)$ we construct a continuous semiflow of continuously differentiable solution operators $x_0\mapsto x_t$, $t\ge0$, on open subsets of the Fr\'echet space $C((-\infty,0],\R^n)$. For nonautonomous equations this yields a continuous process of differentiable solution operators. As an application we obtain processes which incorporate all solutions of Volterra integro-differential equations $x'(t)=\int_0^tk(t,s)h(x(s))ds$. 

\bigskip

\noindent
MSC 2010: 34 K 05, 37 L 05, 45 D 99

\medskip

\noindent
Keywords:  Delay differential equation, unbounded delay, process, Volterra integro-differential equation

\end{abstract}

\maketitle










\section{Introduction}

In the present note we consider the initial value problem
\begin{eqnarray}
x'(t) & = & f(x_t),\\
x_0 & = & \phi\in U,
\end{eqnarray}
for a continuously differentiable map $f:U\to\R^n$ on an open subset $U$ of the Fr\'echet space 
$C=C((-\infty,0],\R^n)$  of continuous maps $(-\infty,0]\to\R^n$, with the topology of locally uniform convergence.
A solution of Eq. (1.1) on an interval $I\subset\R$ is a continuous map $x:(-\infty,0]+I\to\R^n$ so that all segments $x_t:(-\infty,0]\ni s\mapsto x(t+s)\in\R^n$, $t\in  I$, belong to $U$, and 
$x|_I$ is differentiable and satisfies Eq. (1.1) for all $t\in I$. A solution of the initial value problem (IVP)  (1.1)-(1.2)
is a solution on some interval $I=[0,t_x)$, $0<t_x\le\infty$, which satisfies $x_0=\phi$.
Eq. (1.1) generalizes the familiar autonomous delay differential equations, or retarded functional differential equations
\cite{HVL,DvGVLW}, where $U$ is  a subset of a Banach space $C([-r,0],\R^n)$, $r>0$, and covers examples with unbounded delay, including cases of variable, state-dependent delay. In Part I (Sections 2-5) below we show that the IVP (1.1)-(1.2) is well-posed
and that the maximal solutions $x=x^{\phi}$ define a continuous semiflow $\Sigma$  on $U$, by the equation 
$$
\Sigma(t,\phi)=x^{\phi}_t,
$$
with all solution operators $\Sigma(t,\cdot)$ continuously differentiable and their derivatives given by solutions of variational equations. 

Part II (Sections 6-7) deals with nonautonomous equations
\begin{equation}
x'(t)=g(t,x_t)
\end{equation}
for $g:\R\times C\supset V\to\R^n$ continuously differentiable. In Section 6 the result of Part I yields a continuous process of continuously differentiable solution operators.

Among the applications are Volterra integro-differential equations (VIDEs)
\begin{equation}
x'(t)=\int_0^tk(t,s)h(x(s))ds,\quad t>0,
\end{equation}
with $k:\R^2\to\R^{n\times n}$ and $h:\R\to\R^n$ continuously differentiable. Eq. (1.4) can be interpreted as a nonautonomous differential equation with unbounded maximal time-dependent delay $d(t)=t$ at time $t>0$ since $x'(t)$ depends on the values of $x$ for $t-t=0<s<t$  \cite{K0}. In Section 7 we find a  continuously differentiable map $g:\R\times C\to\R^n$ so that solutions of Eq. (1.4) also satisfy Eq. (1.3), which in turn yields a process incorporating all solutions of the VIDE.

The construction of the semiflow $\Sigma$ associated with Eq. (1.1) is a simplified version of a construction in \cite{W7}. It proceeds in the familiar way, via an integral equation for solutions of the IVP (1.1)-(1.2) with the initial data as parameter. However, with the Fr\'echet space $C$ as state space some care has to be taken. This begins with the notion of continuous differentiability. We obtain our result in 2 variants, namely, in the setting of continuous differentiability (1) in the sense of Michal and Bastiani and (2) in the sense of Fr\'echet.  
Let us briefly speak of $C^1_{MB}$-smoothness in case (1) and of $C^1_F$-smoothness in case (2).  
For a continuous map $f:V\supset U\to W$, $V$ and $W$ topological vector spaces and $U\subset V$ open, $C^1_{MB}$-smoothness means that all directional derivatives 
$$
Df(u)v=\lim_{0\neq t\to0}\frac{1}{t}(f(u+tv)-f(u))
$$
exist and that the map
$$
U\times V\ni(u,v)\mapsto Df(u)v\in W
$$
is continuous. By $C^1_F$-smoothness we mean that all directional derivatives exist, that each map $Df(u):V\to W$, $u\in U$, is linear and continuous, and that the map $Df:U\ni u\mapsto Df(u)\in L_c(V,W)$ is continuous with respect to the topology $\beta$ of uniform convergence on bounded sets, on the vector space $L_c(V,W)$ of continuous linear maps $V\to W$.  

In case the topological vector spaces are Banach spaces  $C^1_F$-smoothness is equivalent to the familiar continuous differentiability based on Fr\'echet derivatives, and for finite dimensional spaces $C^1_F$-smoothness and $C^1_{MB}$-smoothness are, of course,  equivalent. In general $C^1_F$-smoothness is the stronger property. \cite{W10} and \cite[Section 8]{W9} contain examples of maps which are $C^1_{MB}$-smooth but not $C^1_F$-smooth. 

The motivation to present results in both settings is that in work involving calculus in topological vector spaces $C^1_{MB}$-smoothness seems to be quite common whereas in our application to VIDEs we obtain an associated equation 
(1.3) with a map $g$ which is in fact $C^1_F$-smooth.

Autonomous equations of the form (1.1) which result from VIDES (1.4) as above, via  an equation of the form (1.3),  are differential equations with unbounded state-dependent delay. As such they are particularly well-behaved, much in contrast to examples with a discrete delay like
$$
x'(t)=F(x(t-d)),\quad d=d(x(t)),
$$
given by continuously differentiable functions $F$ and $d:\R\to(0,\infty)$. For the latter continuously differentiable solution operators exist on submanifolds of the Fr\'echet space
$C^1=C^1((-\infty,0],\R)$  of continuously differentiable maps $(-\infty,0]\to\R$, with the topology of locally uniform convergence of maps and their derivatives \cite{W7,W9}.

For calculus based on $C^1_{MB}$-smoothness we refer to \cite[I.1-I.4]{H}. Below in the appendix Section 8 we collect simple additional facts from calculus based on $C^1_F$-smoothness. Proofs are given in \cite{W9}. 

We have to warn the reader that the hypotheses on continuous differentiability are restrictive in a perhaps surprising way: $C^1_{MB}$-smoothness of a map $f:C\supset U\to\R^n$  implies that $f$ is of {\it locally bounded delay} in the following sense.\\

(lbd) {\it For every $\phi\in U$ there are a neighbourhood $N\subset U$ of $\phi$  and $d>0$ such that for all $\chi,\psi$ in $N$ with
$$
\chi(t)=\psi(t)\quad\text{for all}\quad t\in[-d,0]
$$
we have $f(\chi)=f(\psi)$.}\\

This can be proved in the same way as \cite[Proposition 1.1]{W7}.    

Let us also mention an obvious advantage of the Fr\'echet space $C$ over Banach spaces  of continuous functions
$(-\infty,0]\to\R^n$ which have been used as state spaces \cite{Sch,HMN,W3,MMNN} -  the space $C$ does not exclude segments of solutions by growth or integrability conditions at $-\infty$. Recall that linear autonomous differential equations with constant delay in general have many solutions with arbitrarily fast exponential growth at $-\infty$.

For other work on delay differential equations with solution segments in Fr\'echet spaces of maps $(-\infty,0]\to\R^n$, see \cite{Sch,S}.

{\bf Notation, preliminaries.} 
$\R^{n\times n}$ denotes the vector space of $n\times n$-matrices with real entries.
For basic facts about  topological vector spaces see \cite{R}. Products of topological vector spaces are always equipped with the product topology. We need the following statement on uniform continuity.

\begin{prop}  
\cite[Proposition 1.2]{W9} Suppose $T$ is a topological space, $W$ is a topological vector space, $M$ is a metric space with metric $d$, $g:T\times M\supset U\to W$ is continuous, $U\supset\{t\}\times K$, $K\subset M$ compact. Then $g$ is uniformly continuous on $\{t\}\times K$ in the following sense: For every neighbourhood $N$ of $0$ in $W$ there exist a neighbourhood $T_N$ of $t$ in $T$ and $\epsilon>0$ such that for all $t'\in T_N$, all $\hat{t}\in T_N$,  all $k\in K$, and all $m\in M$ with
$$
d(m,k)<\epsilon\quad\text{and}\quad (t',k)\in U,\quad (\hat{t},m)\in U
$$
we have
$$
g(t',k)-g(\hat{t},m)\in N.
$$
\end{prop}

The vector space of continuous linear maps $V\to W$ between topological vector spaces is denoted by $L_c(V,W)$. The sets
$$
U_{N,B}=\{A\in L_c(V,W):AB\subset N\},
$$
$N$ a neighbourhood of $0$ in $W$ and $B\subset V$ bounded, form a neighbourhood base at $0\in L_c(V,W)$, for the 
topology $\beta$ of uniform convergence on bounded sets.

A Fr\'echet space $F$ is a locally convex topological vector space which is complete and metrizable. The topology is given by a sequence of seminorms 
$|\cdot|_j$, $j\in\N$, which are separating in the sense that $|v|_j=0$ for all $j\in\N$ implies $v=0$. The sets
$$
N_{j,k}=\left\{v\in F:|v|_j<\frac{1}{k}\right\},\quad j\in\N\quad\text{and}\quad k\in\N,
$$
form a neighbourhood base at the origin. If the sequence of seminorms is increasing then the sets
$$
N_j=\left\{v\in F:|v|_j<\frac{1}{j}\right\},\quad j\in\N,
$$
form a neighbourhood base at the origin.

Products of Fr\'echet spaces, closed subspaces of Fr\'echet spaces, and Banach spaces are Fr\'echet spaces.

For a curve, a continuous map $c$ from an interval $I\subset\R$ of positive length into a Fr\'echet space $F$, the
tangent vector at $t\in I$ is 
$$
c'(t)=\lim_{0\neq h\to0}\frac{1}{h}(c(t+h)-c(t))
$$
provided the limit exists. As in \cite[Part I]{H} the curve is said to be continuously differentiable if it has tangent vectors everywhere and if the map
$$
c':I\ni t\mapsto c'(t)\in F
$$
is continuous.

For a continuous map $f:V\supset U\to F$, $V$ and $F$ Fr\'echet spaces and $U\subset V$ open, and for $u\in U,v\in V$ the directional derivative is defined by
$$
Df(u)v=\lim_{0\neq h\to0}\frac{1}{h}(f(u+hv)-f(u))
$$
provided the limit exists. If for $u\in U$ all directional derivatives $Df(u)v$, $v\in V$ exist then the map $Df(u):V\ni v\mapsto Df(u)v\in F$ is called the derivative of $f$ at $u$. 

For continuous maps $f:U\to F$, $V,W,F$ Fr\'echet spaces and $U\subset V\times W$ open,  partial derivatives are defined in the usual way. For example, $D_1f(v,w):V\to F$ is given by
$$
D_1f(v,w)\hat{v}=\lim_{0\neq h\to0}\frac{1}{h}(f(v+h\hat{v},w)-f(v,w)).
$$



The following Fr\'echet spaces are used in the sequel: For $n\in\N$ and $T\ge0$, $C_T=C((-\infty,T],\R^n)$ denotes the Fr\'echet space of continuous maps $(-\infty,T]\to\R^n$ with the seminorms given by
$$
|\phi|_{T,j}=\max_{T-j\le t\le T}|\phi(t)|,\quad\phi\in C_T\quad\text{and}\quad j\in\N,
$$ 
which define the topology of locally uniform convergence. Analogously we consider the space $C_{\infty}=C(\R,\R^n)$, with 
$$
|\phi|_{\infty,j}=\max_{-j\le t\le j}|\phi(t)|.
$$ 
In case $T=0$ we abbreviate $C=C_0$, $|\cdot|_j=|\cdot|_{0,j}$.  In Section 7 on VIDEs we need the Fr\'echet space 
$C^1_{\infty}$ of continuously differentiable maps $\R\to\R^n$, with the seminorms given by $|\phi|_{\infty,1,j}=|\phi|_{\infty,j}+|\phi'|_{\infty,j}$. $C^1$ is the analogous space of continuously differentiable maps $(-\infty,0]\to\R^n$.

The following Banach spaces occur in the sequel:  For $n\in\N$ and $T>0$, $C_{0T}$ denotes the Banach space of continuous maps $[0,T]\to\R^n$ with the norm given by
$$
|\phi|=\max_{0\le t\le T}|\phi(t)|,
$$ 
and $C_{0T,0}$ is the closed subspace of all $\phi\in C_{0T}$ which satisfy $\phi(0)=0$.

The evaluation maps 
$$
E_T:C_T\times(-\infty,T]\to C\quad\text{and}\quad E_{\infty}:C_{\infty}\times\R\to C
$$
given by $(\phi,t)\mapsto\phi_t$ are continuous (see \cite[Proposition 3.1]{W7}), and linear with respect to the first variable. The evaluation map
$$
Ev_{\infty,1}:C^1_{\infty}\times\R\to\R^n,\quad Ev_{\infty,1}(\phi,t)=\phi(t),
$$ 
is $C^1_F$-smooth with 
$$
DEv_{\infty,1}(\phi,t)(\hat{\phi},t_{\ast})=\hat{\phi}(t)+t_{\ast}\phi'(t),
$$
because it  is the composition of the map $E^{10}_{\infty}$ from Proposition 8.8 (\cite[Proposition 9.1 (iii)]{W9}), which is $C^1_F$-smooth, with the evaluation $C\ni\phi\mapsto\phi(0)\in\R^n$, which is linear and continuous The formula for the derivative also follows by means of Proposition 8.8 (\cite[Proposition 9.1 (iii)]{W9}).

For $0\le S<T\le\infty$  the prolongation maps $P_{ST}:C_S\to C_T$ given by 
$(P_{ST}\phi)(t)=\phi(t)$ for $t\le S$ and by $(P_{ST}\phi)(t)=\phi(S)$ for $t>S$ 
are linear and continuous. The same holds for $Z_T:C_{0T,0}\to C_T$ given by
$$
(Z_T\phi)(t)=\phi(t)\quad\text{for}\quad 0\le t\le T,\qquad (Z_T\phi)(t)=0\quad\text{for}\quad t\le 0.
$$

The maps
$$
I_T:C_{0T}\to C_{0T,0},\quad (I_T\phi)(t)=\int_0^t\phi(s)ds,
$$
and
$$
J_T:C_{0T,0}\times C\ni(\chi,\phi)\mapsto P_{0T}\phi+Z_T\chi\in C_T
$$
are linear and continuous.

We reformulate the IVP (1.1)-(1.2) as a fixed point problem in a familiar way: 
Suppose $x$ is a solution of Eq. (1.1) on $[0,T]$ for some $T>0$, with $x_0=\phi\in U$. Then $[0,T]\ni s\to x_s\in C$ is continuous (use $x_s=E_T(x,s)$), and
$$
x(t)-\phi(0)=\int_0^tf(x_s)ds\quad\text{for all}\quad t\in[0,T].
$$
Define $\eta\in C_{0T,0}$ by $\eta(t)=x(t)-\phi(0)$. Then
$$
x_{|(-\infty,T]}=Z_T\eta+P_{0T}\phi,
$$
and
\begin{equation}
\eta(t)=\int_0^tf((Z_T\eta)_s+(P_{0T}\phi)_s)ds\quad\text{for}\quad 0\le t\le T
\end{equation}
which is a fixed point equation for $\eta\in C_{0T,0}$ with parameter $\phi\in U\subset C$.

\bigskip

\begin{center}
PART I
\end{center}

\bigskip

In the following Sections 2-5 we consider an open subset $U\subset C$ and a map $f:U\to\R^n$ which is $C^1_{\ast}$-smooth, $\ast=MB$ or $\ast=F$. 

\bigskip

\section{A substitution operator}

Let $T>0$ in this section. Set
$$
dom_T=\{\xi\in C_T:\xi_t\in U\quad\text{for all}\quad t\in[0,T]\}
$$
and let $F_T:C_T\supset dom_T\to C_{0T}$ be given by
$$
F_T(\xi)(t)=f(\xi_t)\qquad (=f(E_T(\xi,t)))
$$

\begin{prop}
$dom_T$ is open and $F_T$ is continuous.
\end{prop}

\begin{proof}
1. (Openness) Let $\phi\in dom_T$. Due to the continuity of $E_T$ for each $t\in[0,T]$ there are open neighbourhoods
$N_t$ of $\phi$ in $C_T$ and $V_t$ of $t$ in $\R$ with $\psi_s=E_T(\psi,s)\in U$ for all $\psi\in N_t$, $s\in V_t\cap[0,T]$. Due to compactness there exists a finite subset $\tau\subset[0,T]$ with $[0,T]\subset\cup_{t\in\tau}V_t$.
Then $\cap_{t\in\tau}N_t$ is a neighbourhood of $\phi$ in $dom_T$.

\medskip

2. (Continuity) Let $\phi\in dom_T$ and $\epsilon>0$ be given. Apply Proposition 1.1 to the continuous map
$$
dom_T\times[0,T]\ni(\psi,t)\mapsto f(E_T(\psi,t))\in\R^n
$$
and to the compact set $\{\phi\}\times [0,T]$. It follows that there is a neighbourhood $V$ of $\phi$ in $dom_T$ such that for all $\psi\in V$ and for all $t\in[0,T]$ we have
$$
\epsilon>|f(E_T(\psi,t))-f(E_T(\phi,t))|,
$$
hence $\epsilon>|F_T(\psi)-F_T(\phi)|$.
\end{proof}

As $J_T$ is continuous we infer that the set
$$
{\mathcal O}_T=\{(\eta,\phi)\in C_{0T,0}\times C:J_T(\eta,\phi)\in dom_T\}
$$
is open. The fixed point equation (1.5) reads 
\begin{equation}
\eta=(I_T\circ F_T)(J_T(\eta,\phi)),
\end{equation}
for $(\eta,\phi)\in{\mathcal O}_T$. 

\begin{prop}
$F_T$ is $C^1_{\ast}$-smooth, with $(DF_T(\phi)\chi)(t)=Df(\phi_t)\chi_t$.
\end{prop}

\begin{proof}
1. The case $\ast=MB$.

\medskip

1.1. Define 
$$
\Delta:dom_T\times C_T\to C_{0T}.
$$
by $\Delta(\phi,\chi)(t)=Df(\phi_t)\chi_t$. This makes sense because for all $\phi,\chi$ in $C_T$ the map
$$
[0,T]\ni t\mapsto Df(E_T(\phi,t))E_T(\chi,t)\in\R^n
$$
is continuous, due to the continuity of $E_T$ and to the hypothesis that $f$ is $C^1_{MB}$-smooth.

\medskip

Proof that $\Delta$ is continuous : Let $\phi\in dom_T$ and $\chi\in C_T$ be given. Let $\epsilon>0$. Observe that for all $\psi\in dom_T$ and all $\rho\in C_T$ we have
$$
|\Delta_T(\psi,\rho)-\Delta(\phi,\chi)|=\max_{0\le t\le T}|Df(\psi_t)\rho_t-Df(\phi_t)\chi_t)|.
$$
The map
$$
dom_T\times C_T\times[0,T]\ni(\psi,\rho,t)\mapsto Df(\psi_t))\rho_t\in\R^n
$$
is continuous (see the remarks above), hence uniformly continuous on the compact set $\{\phi\}\times\{\chi\}\times[0,T]$.
There is a neighbourhood $N_{\epsilon}$ of $(\phi,\chi)$ in $dom_T\times C_T$ such that for all $(\psi,\rho)\in N_{\epsilon}$ and for all $t\in[0,T]$,
$$
|Df(\psi_t)\rho_t-Df(\phi_t)\chi_t|<\epsilon.
$$
It follows that for all $(\psi,\rho)\in N_{\epsilon}$,
$$
|\Delta(\psi,\rho)-\Delta(\phi,\chi)|\le\epsilon.
$$

1.2. (Directional derivatives) Let $\phi\in dom_T$, $\chi\in C_T$ be given. Choose $r>0$ with $\phi+[-r,r]\chi\in dom_T$. For $0<|h|<r$ , 
$$
\left|\frac{1}{h}(F_T(\phi+h\chi)-F_T(\phi)) - \Delta(\phi,\chi)\right|
$$
$$
=\max_{0\le t\le T}\left|\frac{1}{h}(f(\phi_t+h\chi_t)-f(\phi_t))-Df(\phi_t)\chi_t\right|
$$
$$
=\max_{0\le t\le T}\left|\frac{1}{h}\int_0^1Df(\phi_t+\theta h\chi_t)h\chi_td\theta- Df(\phi_t)\chi_t\right|
$$
$$
=\max_{0\le t\le T}\left|\int_0^1\left[Df(\phi_t+\theta h\chi_t)- Df(\phi_t)\right]\chi_td\theta\right|.
$$
The map
$$
[0,T]\times(-r,r)\times[0,1]\ni(t,h,\theta)\mapsto Df(\phi_t+\theta h\chi_t)\chi_t\in\R^n
$$
is continuous (use $Df(\phi_t+\theta h\chi_t)\chi_t=Df(E_T(\phi+\theta h\chi,t))E_T(\chi,t)$ and the continuity of $E_T$ and the hypothesis that $f$ is $C^1_{MB}$-smooth), hence uniformly continuous on the compact set $[0,T]\times\{0\}\times[0,1]$. Let $\epsilon>0$. Then there exists $\delta_{\epsilon}\in(0,r)$ such that for all $t\in[0,T]$, $h\in(-\delta_{\epsilon},\delta_{\epsilon})$, $\theta\in[0,1]$, we have
\begin{eqnarray*}
\epsilon & > & |Df(\phi_t+\theta h\chi_t)\chi_t- Df(\phi_t+\theta\cdot 0\cdot\chi_t)\chi_t|\\
& = & |Df(\phi_t+\theta h\chi_t)\chi_t- Df(\phi_t)\chi_t|.
\end{eqnarray*}
It follows that for $0<|h|<\delta_{\epsilon}$,
$$
\left|\frac{1}{h}(F_T(\phi+h\chi)-F_T(\phi))-\Delta(\phi,\chi)\right|<\epsilon.
$$
Therefore $DF_T(\phi)\chi$ exists and is equal to $\Delta(\phi,\chi)$. Using Part 1.1 one finds that $F_T$ is $C^1_{MB}$-smooth.

\medskip

2. The case $\ast=F$. Then $f$ is $C^1_{MB}$-smooth, see Proposition 8.2 (\cite[Proposition 3.2]{W9}), and Part 1 above yields that $F_T$ is $C^1_{MB}$-smooth, too. Again by  Proposition 8.2 (\cite[Proposition 3.2]{W9}) it remains to show that the map
$$
C_T\supset dom_T\ni\phi\mapsto DF_T(\phi)\in L_c(C_T,C_{0T})
$$
is continuous with respect to the topology $\beta$ of uniform convergence on bounded subsets of $C_T$. 
Remark 8.1 says that in order to achieve this we have to  
do the following: Given $\xi\in dom_T$, a neighbourhood $V$ of $0$ in $C_{0T}$, and a bounded subset $B\subset C_T$, we have to find a neighbourhood $N$ of $\xi$ in $dom_T$ so that for all $\tilde{\xi}\in N$ and for all $\hat{\xi}\in B$,
$$
[DF_T(\tilde{\xi})-DF_T(\xi)]\hat{\xi}\in V.
$$
We may assume $V=\{\phi\in C_{0T}:|\phi|<\delta\}$ for some $\delta>0$. Then the previous relation follows from
\begin{equation}
\delta>|\{[DF_T(\tilde{\xi})-DF_T(\xi)]\hat{\xi}\}(t)|=|Df(\tilde{\xi}_t)\hat{\xi}_t-Df(\xi_t)\hat{\xi}_t|.
\end{equation}
for all $\tilde{\xi}\in N$, all $\hat{\xi}\in B$, and all $t\in[0,T]$.

\medskip

2.1. So let $\xi\in dom_T$, a bounded set $B\subset C_T$, and $\delta>0$ be given. Proof that 
$$
B_C=\{E_T(\hat{\xi},t)\in C:\hat{\xi}\in B,0\le t\le T\}
$$
is bounded : Let $j\in\N$. We have to show that the seminorm $|\cdot|_j$ is bounded on $B_C$. Choose an integer $k\ge j+T$. The seminorm $|\cdot|_{T,k}$ on $C_T$ is bounded on $B$. For every $\hat{\xi}\in B$ and every $t\in[0,T]$ we infer from
$$
|E_T(\hat{\xi},t)|_j=\max_{-j\le s\le0}|\hat{\xi}(t+s)|\le\max_{-j\le w\le T}|\hat{\xi}(w)|\le|\hat{\xi}|_{T,k}
$$
that $|\cdot|_j$ is bounded on $B_C$.

\medskip

2.2. For every $\tilde{\xi}\in dom_T$,  $\hat{\xi}\in B$, and $t\in[0,T]$ we have
$$
\{Df(\tilde{\xi}_t)-Df(\xi_t)\}\hat{\xi}_t=\{Df(E_T(\tilde{\xi},t))-Df(E_T(\xi,t))\}E_T(\hat{\xi},t),
$$
with $E_T(\hat{\xi},t)\in B_C$. As $f$ is $C^1_F$-smooth and as $E_T$ is continuous the composition
$$
Q:C_T\times\R\supset dom_T\times[0,T]\ni(\tilde{\xi},t)\mapsto Df(E_T(\tilde{\xi},t))\in L_c(C,\R^n)
$$
is continuous with respect to the topology $\beta$ on $L_c(C,\R^n)$. Let $W=\{x\in\R^n:|x|<\delta\}$. The set 
$$
U_{W,B_C}=\{A\in L_c(C,\R^n):AB_C\subset W\}
$$ 
is a neighbourhood of $0$ in $L_c(C,\R^n)$ with respect to the topology $\beta$. Apply Proposition 1.1 (\cite[Proposition 1.2]{W9}) to the map $Q$ and to the compact  set $\{\xi\}\times[0,T]$. It follows that there is a neighbourhood $N$ of $\xi$ in $dom_T\subset C_T$ so that for every $\tilde{\xi}\in N$ and for all $t\in[0,T]
$ the difference
$$
Q(\tilde{\xi},t)-Q(\xi,t)=Df(E_T(\tilde{\xi},t))-Df(E_T(\xi,t))
$$
is contained in $U_{W,B_C}$. Or,
\begin{equation}
\R^n\supset W\ni\{Df(E_T(\tilde{\xi},t))-Df(E_T(\xi,t))\}\hat{\beta}=\{Df(\tilde{\xi}_t)-Df(\xi_t)\}\hat{\beta}_t
\end{equation}
for all $\tilde{\xi}\in N$,  all $t\in[0,T[]$, and all $\hat{\beta}\in B_C\subset C$. For every  $\tilde{\xi}\in N$, $t\in[0,T[]$, and $\hat{\xi}\in B$ we have $(\hat{\beta}=)\quad\hat{\xi}_t\in B_C$. Using the relation (2.3) we finally obtain
the inequality (2.2).
\end{proof}

It follows that the map $B_T:{\mathcal O}_T\to C_{0T,0}$ given by
$$
B_T(\eta,\phi)=(I_T\circ F_T)(J_T(\eta,\phi))
$$
is $C^1_{\ast}$-smooth.

\section{Uniform contractions and local solutions}

In order to prepare the proof that certain restrictions of $B_T$, for $T>0$ sufficiently small, are uniform contractions observe first that for $T>0$ and  $(\eta,\phi)$ and $(\hat{\eta},\phi)$ both in ${\mathcal O_T}$ we have
\begin{eqnarray*}
|B_T(\hat{\eta},\phi)-B_T(\eta,\phi)| & = & |I_T(F_T(J_T(\hat{\eta},\phi)))-I_T(F_T(J_T(\eta,\phi)))|\\
& = & |I_T\{F_T(J_T(\hat{\eta},\phi))-F_T(J_T(\eta,\phi))\}|\\
& \le & T\,\max_{0\le t\le T}|\{F_T(J_T(\hat{\eta},\phi))-F_T(J_T(\eta,\phi))\}(t)|,
\end{eqnarray*}
and for all $t\in[0,T]$
$$
\{F_T(J_T(\hat{\eta},\phi))-F_T(J_T(\eta,\phi))\}(t)=f((P_{0T}\phi)_t+(Z_T\hat{\eta})_t)-f((P_{0T}\phi)_t+(Z_T\eta)_t).
$$
In case the line segment between the arguments of $f$ belongs to $U$ the last term equals
$$
\int_0^1Df((P_{0T}\phi)_t+(Z_T\eta)_t+\theta[(Z_T\hat{\eta})_t-(Z_T\eta)_t])[(Z_T\hat{\eta})_t-(Z_T\eta)_t]d\theta.
$$

\begin{prop}
Let $\phi\in dom_T$ be given. There exist $T=T_{\phi}>0$, a neighbourhood $V=V_{\phi}$ of $\phi$ in $dom_T$, $\epsilon=\epsilon_{\phi}>0$, and $j=j_{\phi}\in\N$ such that for all $S\in(0,T]$, all $\chi\in V$, all $\eta$ and $\tilde{\eta}$ in $C_{0S,0}$ with $|\eta|<\epsilon$ and $|\tilde{\eta}|<\epsilon$, all $w\in[0,S]$, and all $\theta\in[0,1]$, we have
\begin{equation}
(P_{0S}\chi)_w+(Z_S\eta)_w+\theta[(Z_S\tilde{\eta})_w-(Z_S\eta)_w]\in U
\end{equation}
and
$$
|Df((P_{0S}\chi)_w+(Z_S\eta)_w+\theta[(Z_S\tilde{\eta})_w-(Z_S\eta)_w])[(Z_S\tilde{\eta})_w-(Z_S\eta)_w)]
\le 2j\,|\tilde{\eta}-\eta|.
$$
\end{prop}

\begin{proof}
1. Let $\phi\in U$ be given. As $f$ is $C^1_{MB}$-smooth the map $U\times C\ni(\chi,\eta)\mapsto Df(\chi)\eta\in\R^n$ is continuous. Then there are
neighbourhoods $V'$ of $\phi$ in $U$ and $N$ of $0$ in $C$ with
$$
|Df(\chi)\eta|=|Df(\chi)\eta-Df(\phi)0|<1\quad\text{for all}\quad \chi\in V',\eta\in N.
$$
There exists $j=j_N\in\N$ with
$$
\left\{\zeta\in C:|\zeta|_j<\frac{1}{j}\right\}\subset N.
$$

2. By the continuity of the map 
$$
\R\ni t\mapsto E_{\infty}(P_{0\infty}\phi,t)\in C
$$ 
at $t=0$, with $E_{\infty}(P_{0\infty}\phi,0)=\phi$, there exists $T>0$ with $E_{\infty}(P_{0\infty}\phi,t)\in V'$ for all $t\in[0,T]$. The continuous map
$$
\alpha:C\times C_{0T,0}\times[0,T]\ni(\chi,\eta,t)\mapsto E_{\infty}(P_{0\infty}\chi,t)+E_T(Z_T\eta,t)\in C
$$
satisfies $\alpha(\phi,0,t)=E_{\infty}(P_{0\infty}\phi,t)\in V'$ for all $t\in[0,T]$ and is uniformly continuous on the compact set $\{\phi\}\times\{0\}\times[0,T]$. It follows that there exist a neighbourhood $V$ of $\phi$ in $V'$ and $\epsilon>0$ such that
$$
E_{\infty}(P_{0\infty}\chi,t)+E_T(Z_T\eta,t)=\alpha(\chi,\eta,t)\in V'
$$
for all $\chi\in V$, $\eta\in C_{0T,0}$ with $|\eta|<\epsilon$, and $t\in[0,T]$.
Observe that $E_{\infty}(P_{0\infty}\chi,t)=E_T(P_{0T}\chi,t)$ for these $\chi$ and $t$. 

\medskip

3. Let $0<S<T$ and let $\chi\in V$, $\eta\neq\tilde{\eta}$ in $C_{0S,0}$ be given, with $|\eta|<\epsilon$ and  
$|\tilde{\eta}|<\epsilon$. Let $0\le w\le S$, $0\le\theta\le1$. Then
$$
|P_{ST}\eta|\le|\eta|<\epsilon\quad\text{and}\quad|P_{ST}\tilde{\eta}|\le|\tilde{\eta}|<\epsilon.
$$
By convexity,
$$
|P_{ST}\eta+\theta[P_{ST}\tilde{\eta}-P_{ST}\eta]|<\epsilon.
$$
The choice of $V$ and $\epsilon$ in Part 2 yields
$$
V'\ni E_{\infty}(P_{0\infty}\chi,w)+E_T(Z_T(P_{ST}\eta+\theta[P_{ST}\tilde{\eta}-P_{ST}\eta]),w).
$$
Due to $0\le w\le S$,
$$
E_T(Z_TP_{ST}\eta,w)=(Z_S\eta)_w\quad\text{and}\quad E_T(Z_TP_{ST}\tilde{\eta},w)=(Z_S\tilde{\eta})_w
$$
and
$$
E_T(Z_T(P_{ST}\eta+\theta[P_{ST}\tilde{\eta}-P_{ST}\eta]),w)
$$
$$
= E_T(Z_TP_{ST}\eta,w)+\theta[E_T(Z_TP_{ST}\tilde{\eta},w)-E_T(Z_TP_{ST}\eta,w)]
$$
$$
=(Z_S\eta)_w+\theta[(Z_S\tilde{\eta})_w-(Z_S\eta)_w].
$$
Using this and $E_{\infty}(P_{0\infty}\chi,w)=(P_{0S}\chi)_w$ one arrives at
$$
U\supset V'\ni (P_{0S}\chi)_w+(Z_S\eta)_w+\theta[(Z_S\tilde{\eta})_w-(Z_S\eta)_w].
$$

4. For
$$
\zeta=\frac{1}{2j|\eta-\tilde{\eta}|}(\tilde{\eta}-\eta)\in C_{0S,0}.
$$
we have
\begin{eqnarray*}
|(Z_S\zeta)_w|_j & = & \max_{-j\le t\le0}|(Z_S\zeta)(w+t)|=\max_{w-j\le s\le w}|(Z_S\zeta)(s)|\\
& \le & \max_{0\le s\le S}|(Z_S\zeta)(s)|=\max_{0\le s\le S}|\zeta(s)|=|\zeta|<\frac{1}{j},
\end{eqnarray*}
hence $(Z_S\zeta)_w\in N$. Using this and the result of Part 3 we infer
$$
1>|Df((P_{0S}\chi)_w+(Z_S\eta)_w+\theta[(Z_S\tilde{\eta})_w-(Z_S\eta)_w])(Z_S\zeta)_w|
$$
$$
= |Df((P_{0S}\chi)_w+(Z_S\eta)_w+\theta[(Z_S\tilde{\eta})_w-(Z_S\eta)_w])\frac{1}{2j|\eta-\tilde{\eta}|}
(Z_S(\tilde{\eta}-\eta))_w|
$$
$$
=|Df((P_{0S}\chi)_w+(Z_S\eta)_w+\theta[(Z_S\tilde{\eta})_w-(Z_S\eta)_w])\frac{1}{2j|\eta-\tilde{\eta}|}
((Z_S\tilde{\eta})_w-(Z_S\eta)_w)|
$$
which implies the estimate in the proposition.
\end{proof}

\medskip

Let $\phi\in U$, and let $T=T_{\phi}>0$, a convex neighbourhood $V=V_{\phi}$ of $\phi$ in $U$, $\epsilon=\epsilon_{\phi}>0$, and $j=j_{\phi}\in\N$ be given as in Proposition  3.1.

\begin{prop}
For every $S\in(0,T)$, $\chi\in V$, $\eta$ and $\tilde{\eta}$ in $C_{0S,0}$ with $|\eta|<\epsilon$ and
$|\tilde{\eta}|<\epsilon$, we have
$$
(\eta,\chi)\in{\mathcal O}_S\,\,(\tilde{\eta},\chi)\in{\mathcal O}_S,\quad\text{and}\quad
|B_S(\tilde{\eta},\chi)-B_S(\eta,\chi)|\le 2jS|\tilde{\eta}-\eta|.
$$
\end{prop}

\begin{proof}
Let  $S\in(0,T)$, $\chi\in V$, $\eta$ and $\tilde{\eta}$ in $C_{0S,0}$ with $|\eta|<\epsilon$ and
$|\tilde{\eta}|<\epsilon$ be given. The relation (3.1) for $0\le w\le S$, with $\theta=0$ and $\theta=1$, yields 
$(\eta,\chi)\in{\mathcal O}_S$ and $(\tilde{\eta},\chi)\in{\mathcal O}_S$. Also, by the same argument, for every $\theta\in[0,1]$,
\begin{equation}
dom_S\ni P_{0S}\chi+Z_S\eta+\theta[Z_S\tilde{\eta}-Z_S\eta]=J_S(\eta,\chi)+\theta[J_S(\tilde{\eta},\chi)-J_s(\eta,\chi)].
\end{equation}
We have,
\begin{eqnarray*}
|B_S(\tilde{\eta},\chi)-B_S(\eta,\chi)| & = & |I_S[F_S(J_S(\tilde{\eta},\chi))-F_S(J_S(\eta,\chi))]|\\
& \le & S\,\max_{0\le w\le S}|F_S(J_S(\tilde{\eta},\chi))(w)-F_S(J_S(\eta,\chi))(w)|\\
& = & S|F_S(J_S(\tilde{\eta},\chi))-F_S(J_S(\eta,\chi))|.
\end{eqnarray*}
As $F_S$ is $C^1_{MB}$-smooth and as the relation (3.2) holds for all $\theta\in[0,1]$ we obtain 
that the last term equals
$$
S\,\left|\int_0^1DF_S(J_S(\eta,\chi)+\theta[J_S(\tilde{\eta},\chi)-J_S(\eta,\chi)])[J_S(\tilde{\eta},\chi)-J_S(\eta,\chi)]d\theta\right|
$$
$$
=S\left|\int_0^1DF_S(P_{0S}\chi+Z_S\eta+\theta[Z_S\tilde{\eta}-Z_S\eta])[Z_S\tilde{\eta}-Z_S\eta]d\theta\right|
$$
$$
\le S\,\max_{0\le\theta\le1}(\max_{0\le w\le S}|Df((P_{0S}\chi)_w+(Z_S\eta)_w+\theta[(Z_S\tilde{\eta})_w-(Z_S\eta)_w])[(Z_S\tilde{\eta})_w-(Z_S\eta)_w]|)
$$
$$
\le S\cdot 2j\cdot|\tilde{\eta}-\eta|\qquad\text{(by Proposition 3.1)}
$$
\end{proof}

\medskip

\begin{prop}
$\lim_{S\searrow0}B_S(0,\phi)=0$.
\end{prop}

\begin{proof}
Use
\begin{eqnarray*}
|B_S(0,\phi)| & = & |I_S(F_S(J_S(0,\phi)))|\le S|F_S(J_S(0,\phi))|=S\,\max_{0\le w\le S}|f((P_{0S}\phi)_w)|\\
& \le & S\,\max_{0\le w\le T}|f((P_{0T}\phi)_w)|.
\end{eqnarray*}
\end{proof}

\medskip

\begin{prop}
There exist $S_{\phi}\in(0,T_{\phi})$ and an open neighbourhood $W_{\phi}$ of $\phi$ in $V_{\phi}$ such that for all $\chi\in W_{\phi}$, for all $S\in(0,S_{\phi}]$, and all
$\eta\in C_{0S,0}$ and $\tilde{\eta}\in C_{0S,0}$ with $|\eta|\le\frac{\epsilon_{\phi}}{2}$ and  $|\tilde{\eta}|\le\frac{\epsilon_{\phi}}{2}$, we have
$$
(\eta,\chi)\in {\mathcal O}_S,\,\,(\tilde{\eta},\chi)\in {\mathcal O}_S,
$$
$$
|B_S(\eta,\chi)|<\frac{\epsilon_{\phi}}{2}\quad\text{and}\quad
|B_S(\tilde{\eta},\chi)-B_S(\eta,\chi)|\le\frac{1}{2}|\tilde{\eta}-\eta|.
$$
\end{prop}

\begin{proof}
1. Choose $S_{\phi}\in(0,T_{\phi})$ with
$$
|B_S(0,\phi)|<\frac{\epsilon_{\phi}}{8}\quad\text{for all}\quad S\in(0,S_{\phi}],
$$
which is possible due to Proposition 3.3, and
$$
2jS_{\phi}<\frac{1}{2}.
$$
As $B_{S_{\phi}}$ is continuous there exists an open neigbourhood $W_{\phi}$ of $\phi$ in $V_{\phi}$ so that for all $\chi\in W_{\phi}$,
$$
|B_{S_{\phi}}(0,\chi)-B_{S_{\phi}}(0,\phi)|<\frac{\epsilon_{\phi}}{8}.
$$

2. Now let $S\in(0,S_{\phi}]$ be given. For every $\chi\in W_{\phi}$ and $t\in[0,S]$,
$$
B_S(0,\chi)(t)=\int_0^tf((P_{0S}\chi)_w)dw=\int_0^tf((P_{0S_{\phi}}\chi)_w)dw=B_{S_{\phi}}(0,\chi)(t).
$$
Using this (for $\chi$ and $\phi$) one gets
$$
|B_S(0,\chi)-B_S(0,\phi)|\le|B_{S_{\phi}}(0,\chi)-B_{S_{\phi}}(0,\phi)|<\frac{\epsilon_{\phi}}{8}
$$

3. Let $\chi\in W_{\phi}$, $\eta\in C_{0S,0}$, $\tilde{\eta}\in C_{0S,0}$ be given, with $|\eta|\le\frac{\epsilon_{\phi}}{2}$ and  $|\tilde{\eta}|\le\frac{\epsilon_{\phi}}{2}$. Proposition 3.2 yields
$$
|B_S(\tilde{\eta},\chi)-B_S(\eta,\chi)|\le 2jS|\tilde{\eta}-\eta|\le\frac{1}{2}\tilde{\eta}-\eta|.
$$
Furthermore,
\begin{eqnarray*}
|B_S(\eta,\chi)| & \le & |B_S(\eta,\chi)-B_S(0,\chi)|+|B_S(0,\chi)-B_S(0,\phi)|\\
& & +|B_S(0,\phi)|\\
& < & \frac{1}{2}|\eta|+\frac{\epsilon_{\phi}}{8}+\frac{\epsilon_{\phi}}{8}\le\frac{1}{2}\frac{\epsilon_{\phi}}{2}+\frac{2\epsilon_{\phi}}{8}=
\frac{\epsilon_{\phi}}{2}.
\end{eqnarray*}
\end{proof}

\medskip

Let $S\in(0,S_{\phi}]$ be given. In  case $\ast=MB$ the uniform contraction result \cite[Theorem 7.2]{W7} applies to the map 
$$
\{\eta\in C_{0S,0}:|\eta|<\epsilon_{\phi}\}\times W_{\phi}\ni(\eta,\chi)\mapsto B_S(\eta,\chi)\in C_{0S,0},
$$
with $M=M_{\phi}=\{\eta\in C_{0S,0}:|\eta|\le\frac{\epsilon_{\phi}}{2}\}$. In case $\ast=F$ the uniform contraction result Theorem 8.7 (\cite[Theorem 5.2]{W9}) applies to the same map and to the same set $M$. It follows that the relations
$$
B_S(\eta,\chi)=\eta\in M,\quad\chi\in W_{\phi}
$$
define a map
$$
W_{\phi}\ni\chi\mapsto\eta_{\chi}\in C_{0S,0}
$$
which is  $C^1_{\ast}$-smooth. As the maps $P_{0S}$ and $Z_S$ are linear and continuous it follows that the map
$$
\Sigma_{\phi}:W_{\phi}\ni\chi\mapsto P_{0S}\chi+Z_S\eta_{\chi}\in C_S
$$
is $C^1_{\ast}$-smooth. Using this and the continuous linear maps $E_S(\cdot,t):C_S\to C$, $0\le t\le S$, one gets that each map
$$
W_{\phi}\ni\chi\mapsto E_S(\Sigma_{\phi}(\chi),t)\in C,\quad 0\le t\le S,
$$
is $C^1_{\ast}$-smooth. The map
$$
[0,S]\times W_{\phi}\ni(t,\chi)\mapsto E_S(\Sigma_{\phi}(\chi),t)\in C
$$ 
is continuous.

\medskip

\begin{prop}
Let $S\in(0,S_{\phi}]$ and $\chi\in W_{\phi}$ be given. The map $x=x^{(\chi)}=\Sigma_{\phi}(\chi)$ is  a solution of Eq. (1.1) on $[0,S]$, with $x_0=\chi$.
\end{prop}

\begin{proof}
$x=\Sigma_{\phi}(\chi)\in C_S$ is continuous, with 
$$
x_0=\Sigma_{\phi}(\chi)_0=(P_{0S}\chi)_0+(Z_S\eta_{\chi})_0=\chi+0=\chi.
$$
For $0\le t\le S$, 
\begin{eqnarray*}
x(t) & = & (P_{0S}\chi)(t)+(Z_S\eta_{\chi})(t)=\chi(0)+\eta_{\chi}(t)\\
& = & \chi(0)+B_S(\eta_{\chi},\chi)(t)=\chi(0)+\int_0^tf((P_{0S}\chi+Z_S\eta_{\chi})_w)dw\\
& = & \chi(0)+\int_0^tf(E_S(\Sigma_{\phi}(\chi),w))dw.
\end{eqnarray*}
The last integrand is continuous. It follows that the restriction $x|_{[0,S]}$  is continuously differentiable, with
$$
(x|_{[0,S]})'(t)=f((\Sigma_{\phi}(\chi))_t)=f(x_t)\quad\text{for all}\quad t\in[0,S].
$$
\end{proof}

From the remarks preceding Proposition 3.5 we see that all maps 
$$
W_{\phi}\ni\chi\mapsto x^{(\chi)}_t\in C,\quad0\le t\le S,
$$
are $C^1_{\ast}$-smooth, and that the map
$$
[0,S]\times W_{\phi}\ni(t,\chi)\mapsto x^{(\chi)}_t\in C
$$
is continuous.

\begin{prop}
(Uniqueness) Suppose $x$ is a  solution of Eq. (1.1) on the interval $I$ and $\tilde{x}$ is  a  solution of Eq. (1.1) on the interval $\tilde{I}$, both of positive length, and $0=\min\,I=\min\,\tilde{I}, x_0=\tilde{x}_0$. Then $x(t)=\tilde{x}(t)$ on $I\cap\tilde{I}$.
\end{prop}

\begin{proof}
1. Proof that there exists $\tau>0$ with $[0,\tau]\subset I\cap\tilde{I}$ and $x(t)=\tilde{x}(t)$ for all $t\le\tau$.
Let $\phi=x_0\quad(=\tilde{x}_0\in U)$. Consider $T_{\phi},\epsilon_{\phi},S_{\phi}$ as in Proposition 3.4. By continuity there exists $\tau=S\in(0,S_{\phi}]\cap I\cap\tilde{I}$ such that for $0\le t\le S$,
$$
|x(t)-\phi(0)|<\frac{\epsilon_{\phi}}{2}\quad\text{and}\quad|\tilde{x}(t)-\phi(0)|<
\frac{\epsilon_{\phi}}{2}.
$$
Define
\begin{eqnarray*}
y & = & x|_{(-\infty,S]}-P_{0S}\phi,\quad\eta=y|_{[0,S]}\in C_{0S,0},\\
\tilde{y} & = & \tilde{x}|_{(-\infty,S]}-P_{0S}\phi,\quad\tilde{\eta}=\tilde{y}|_{[0,S]}\in C_{0S,0}.
\end{eqnarray*}
Then
$$
|\eta|<\frac{\epsilon_{\phi}}{2}\quad\text{and}\quad|\tilde{\eta}|<\frac{\epsilon_{\phi}}{2},
$$
and for $0\le t\le S$,
\begin{eqnarray*}
B_S(\eta,\phi)(t) & = & \int_0^tf((P_{0S}\phi)_w+(Z_S\eta)_w)dw\\
& = & \int_0^tf(x_w)dw=x(t)-\phi(0)=\eta(t).
\end{eqnarray*}
Hence $B_S(\eta,\phi)=\eta$. Analogously, $B_S(\tilde{\eta},\phi)=\tilde{\eta}$. Proposition 3.4 yields
$$
|\tilde{\eta}-\eta|=|B_S(\tilde{\eta},\phi)-B_S(\eta,(\phi)|\le\frac{1}{2}|\tilde{\eta}-\eta|,
$$
which gives $\tilde{\eta}=\eta$ and thereby $\tilde{x}(t)=x(t)$ on $[0,S]=[0,\tau]$.

2. The interval $J=I\cap\tilde{I}$ has positive length, and $\min\,J=0$. Assume $x(u)\neq\tilde{x}(u)$ for some $u\in J$.
Then $0<u$, and by continuity, $t_J=\inf\{t\in J:x(t)\neq\tilde{x}(t)\}<u\le\sup J$. On $(-\infty,t_J]$ we have $x(t)=\tilde{x}(t)$ while every neighbourhood of $t_J$ contains $t>t_J$ in $J$ with $x(t)\neq\tilde{x}(t)$.
The continuously differentiable function $y:(-\infty,\sup\,J-t_J)\to\R^n$ given by $y(t)=x(t+t_J)$ satisfies
$$
y'(t)=x'(t+t_J)=f(x_{t+t_J})=f(y_t)
$$
for $0\le t<\sup\,J-t_J$ (with the right derivative at $t=0$). Analogously the function $\tilde{y}:(-\infty,\sup\,J-t_J)\to\R^n$
given by $y(t)=\tilde{x}(t+t_J)$ is a solution of Eq. (1.1) on $[0,\sup\,J-t_J)$, and $y_0=\tilde{y}_0$. Part 1 of the proof yields $y(t)=\tilde{y}(t)$ on an interval $[0,\tau]$ with $0<\tau<\sup\,J-t_J$. This implies $x(t)=\tilde{x}(t)$ on $[t_J,t_J+\tau]$, contradicting the definition of $t_J$. 
\end{proof}

\section{The semiflow of continuously differentiable solution operators}

Now we proceed as in \cite[Section 5]{W7}. Proofs are included for convenience.  The maximal solution of the IVP (1.1)-(1.2) given by the initial condition $x_0=\phi\in U$ is defined as follows. Set
$$
t_{\phi}=\sup\{t>0:\text{There is a solution of Eq. (1.1) on}\quad[0,t]\quad\text{with}\quad x_0=\phi\}\le\infty.
$$
By Proposition 3.5, $0<t_{\phi}$. Using Proposition 3.6 one obtains a solution $x^{\phi}$ of Eq. (1.1) on $[0,t_{\phi})$, with
$x^{\phi}_0=\phi$, by
$$
x^{\phi}(t)=x(t)
$$ 
for $0<t<t_{\phi}$, where $x$ is any solution of Eq. (1.1) on $[0,t']$ with $t<t'<t_{\phi}$ and $x_0=\phi$. 

It is easy to show that any solution of Eq. (1.1) on some interval $I$ of positive length with $\min\,I=0$ and $x_0=\phi$ is a restriction of $x^{\phi}$.

Set
$$
\Omega=\{(t,\phi)\in[0,\infty)\times U:t<t_{\phi}\}
$$
and define $\Sigma:\Omega\to U$ by $\Sigma(t,\phi)=x^{\phi}_t$.


\begin{prop}
(Semiflow) $\{0\}\times U\subset\Omega$, $\Sigma(0,\phi)=\phi$ for all $\phi\in U$, and if $(t,\phi)\in\Omega$ and $(s,\Sigma(t,\phi))\in\Omega$ then
$$
(s+t,\phi)\in\Omega\quad\text{and}\quad\Sigma(s,\Sigma(t,\phi))=\Sigma(s+t,\phi).
$$
\end{prop}

\begin{proof}
For every $\phi\in U$, $0<t_{\phi}$, hence $(0,\phi)\in\Omega$ and $\Sigma(0,\phi)=x^{\phi}_0=\phi$. Let $(t,\phi)\in\Omega$ and  $(s,\Sigma(t,\phi))\in\Omega$. Let $x=x^{\phi},\psi=x_t,y=x^{\psi}$. Define
$\xi:(-\infty,s+t]\to\R^n$ by $\xi(u)=y(u-t)$. For $u\le t$ we get
$$
\xi(u)=y(u-t)=\psi(u-t)=x_t(u-t)=x(u).
$$ 
In particular, $\xi_0=\phi$ and $\xi'(u)=f(\xi_u)$ for $0\le u\le t$ (with the right derivative at $u=0$) . For $t<u\le t+s$,
$$
\xi'(u)=y'(u-t)=f(y_{u-t})=f(\xi_u).
$$
It follows that $\xi$ is a restriction of $x^{\phi}$. Hence $s+t<t_{\phi}$, or, $(s+t,\phi)\in\Omega$, and
$$
\Sigma(s+t,\phi)=\xi_{s+t}=y_s=\Sigma(s,\psi)=\Sigma(s,\Sigma(t,\phi)).
$$
\end{proof}

For $t\ge0$ with $\Omega_t=\{\phi\in U:(t,\phi\in\Omega\}\neq\emptyset$ consider the solution operator
$$
\Sigma_t:\Omega_t\to U
$$
given by $\Sigma_t(\phi)=\Sigma(t,\phi)$.

\begin{prop}
For every $(t,\phi)\in\Omega$ there exist an open neighbourhood $N\subset U$ of $\phi$ and $\epsilon>0$ with $[0,t+\epsilon)\times N\subset\Omega$, $\Sigma|_{[0,t+\epsilon)\times N}$ continuous, and $\Sigma_t|_N$
$C^1_{\ast}$-smooth.
\end{prop}

\begin{proof}
1. Let $(t,\phi)\in\Omega$ be given. The remarks following Proposition 3.5 show that $t=0$ is contained in the set
\begin{eqnarray*}
A & = & \{s\in[0,t_{\phi}):\text{There exist an open neighbourhood}\,\,V_s\subset U\,\,\text{of}\,\,\phi\\
& & \text{and}\,\,\epsilon_s>0\,\,\text{with}\,\,[0,s+\epsilon_s)\times V_s\subset\Omega,\,\,\Sigma|_{[0,s+\epsilon_s)\times V_s}\,\,\text{continuous,}\\
& & \text{and}\,\,\Sigma_s|_{V_s}\,\,C^1_{\ast}-\text{smooth}\}.
\end{eqnarray*}
Let $t_A=\sup\,A\le t_{\phi}$. It remains to prove that $t_A=t_{\phi}$. 

2. Suppose $t_A<t_{\phi}$. Set $\psi=\Sigma(t_A,\phi)$. Again by the remarks following Proposition 3.5, there exist an open neighbourhood $W\subset U$ of $\psi$ and $\tau>0$ with $[0,\tau]\times W\subset\Omega$ so that
$\Sigma|_{[0,\tau]\times W}$ is continuous and all $\Sigma_u|_W$, $0\le u\le\tau$, are $C^1_{\ast}$-smooth.
The flowline $[0,t_{\phi})\ni s\mapsto x^{\phi}_s\in U$ is continuous (observe $x^{\phi}_s=E_u(x^{\phi}|_{(-\infty,u]},s)$ for $0\le s<u<t_{\phi}$, with $E_u$ continuous). It follows that there exists 
$$
t_0\in A\cap\left(t_A-\frac{\tau}{2},t_A\right)\quad\text{with}\quad x^{\phi}_{t_0}\in W.
$$
From $t_0\in A$ one obtains an open neighbourhood $N_0\subset U$ of $\phi$ and $\epsilon_0>0$ so that
$[0,t_0+\epsilon_0)\times N_0\subset\Omega$, and $\Sigma|_{[0,t_0+\epsilon_0)\times N_0}$ is continuous, and
$\Sigma_{t_0}|_{N_0}$ is $C^1_{\ast}$-smooth. Because of continuity and $x^{\phi}_{t_0}\in W$ one may assume 
$\Sigma_{t_0}(N_0)\subset W$. For $t_0<u<t_A+\frac{\tau}{2}$ and $\chi\in N_0$,
$$
0<u-t_0<\tau\quad\text{and}\quad\Sigma_{t_0}(\chi)\in W,
$$ 
which gives $(u,\chi)=((u-t_0)+t_0,\chi)\in\Omega$ and
$$
\Sigma(u,\chi)=\Sigma(u-t_0,\Sigma(t_0,\chi)).
$$
It follows that $\Sigma|_{(t_0,t_A+\frac{\tau}{2})\times N_0}$ is continuous, which in combination with the continuity of the restriction $\Sigma|_{[0,t_0+\epsilon_0)\times N_0}$ yields that the restriction of $\Sigma$ to $[0,t_A+\frac{\tau}{2})\times N_0$ is continuous.

3. For $u=t_A+\frac{\tau}{4}$ and $\chi\in N_0$,
$$
\Sigma(u,\chi)=\Sigma(u-t_0,\Sigma(t_0,\chi))=\Sigma_{u-t_0}\circ \Sigma_{t_0}(\chi)
$$
with $0<u-t_0<\tau$. Recall $\Sigma_{t_0}(N_0)\subset W$. Now it follows that $\Sigma_u|_{N_0}$ is $C^1_{\ast}$-smooth. Combining this with the result of Part 2 of the proof one concludes that $u>t_A$ belongs to $A$, contradicting
$t_A=\sup\,A$.
\end{proof}

\begin{cor}
The semiflow $\Sigma$ is continuous, each set $\Omega_t$, $t\ge0$, is open in $X_f$, and each solution operator $\Sigma_t$, $t\ge0$ and $\Omega_t\neq\emptyset$, is $C^1_{\ast}$-smooth.
\end{cor}

\begin{proof} 
Let $t\ge0$ and $\phi\in\Omega_t$ be given. Then $(t,\phi)\in\Omega$, and for $N$ chosen according to Proposition 4.2 we get $N\subset\Omega_t$. This shows that $\Omega_t\subset U$ is an open subset of $C$. The remaining assertions are obvious from Proposition 4.2.
\end{proof}

\section{Linearized solution operators and the variational equation}

For $\phi\in U$ the derivatives $D\Sigma_t(\phi):C\to C$, $0\le t<t_{\phi}$, are given by a variational equation.
The proof requires the following version of \cite[Proposition 5.5]{W7}.

\begin{prop}
Let $\phi\in U$, $0\le t<t_{\phi}$, $\hat{\phi}\in C$, and $s\le0$. Then
\begin{eqnarray*}
(D\Sigma_t(\phi)\hat{\phi})(s) & = & \hat{\phi}(t+s)\quad\text{in case}\quad t+s\le0,\\
(D\Sigma_t(\phi)\hat{\phi})(s) & = & (D\Sigma_{t+s}(\phi)\hat{\phi})(0)\quad\text{in case}\quad0\le t+s.
\end{eqnarray*}
\end{prop}

\begin{proof}
Each linear map 
$$
ev_s:C\ni\psi\mapsto\psi(s)\in\R^n,\quad s\le0,
$$ 
is continuous. Let $\phi\in U$, $0\le t<t_{\phi}$, $\hat{\phi}\in C$, $s\le0$. 
Then
\begin{eqnarray*}
(D\Sigma_t(\phi)\hat{\phi})(s) & = & ev_s(D\Sigma_t(\phi)\hat{\phi})=D(ev_s\circ \Sigma_t)(\phi)\hat{\phi}\\
& = & D\{\Omega_t\ni\tilde{\phi}\mapsto x^{\tilde{\phi}}_t(s)\in\R^n\}(\phi)\hat{\phi}\\
& = & D\{\Omega_t\ni\tilde{\phi}\mapsto x^{\tilde{\phi}}(t+s)\in\R^n\}(\phi)\hat{\phi}.
\end{eqnarray*}
In case $0\le t+s$ the set $\Omega_t\subset\Omega_{t+s}$ is an open neighbourhood of $\phi$ in $U$, and 
\begin{eqnarray*}
D\{\Omega_t\ni\tilde{\phi}\mapsto x^{\tilde{\phi}}(t+s)\in\R^n\}(\phi)\hat{\phi} & = & 
D\{\Omega_t\ni\tilde{\phi}\mapsto x^{\tilde{\phi}}_{t+s}(0)\in\R^n\}(\phi)\hat{\phi}\\
& = & D(ev_0\circ\Sigma_{t+s})(\phi)\hat{\phi}\\
& = & ev_0(D\Sigma_{t+s}(\phi)\hat{\phi})=(D\Sigma_{t+s}(\phi)\hat{\phi})(0)
\end{eqnarray*}
while in case $t+s\le0$,
\begin{eqnarray*}
D\{\Omega_t\ni\tilde{\phi}\mapsto x^{\tilde{\phi}}(t+s)\in\R^n\}(\phi)\hat{\phi} & = & 
D\{\Omega_t\ni\tilde{\phi}\mapsto \tilde{\phi}(t+s)\in\R^n\}(\phi)\hat{\phi}\\
& = & D\,ev_{t+s}(\phi)\hat{\phi}=ev_{t+s}(\hat{\phi})=\hat{\phi}(t+s).
\end{eqnarray*}
\end{proof}

Now we follow \cite[Section 6]{W7}. For $\phi\in U$  define the map $v^{\phi,\hat{\phi}}:(-\infty,t_{\phi})\to\R^n$ by
\begin{eqnarray*}
v^{\phi,\hat{\phi}}(t) & = & (D\Sigma_t(\phi)\hat{\phi})(0)\quad\text{for}\quad0\le t<t_{\phi},\\
v^{\phi,\hat{\phi}}(t) & = & \hat{\phi}(t)\quad\text{for}\quad t<0.
\end{eqnarray*}

\begin{prop}
Let $\phi\in U$ and $\hat{\phi}\in C$ be given and consider the map $v=v^{\phi,\hat{\phi}}$. For every $t\in[0,t_{\phi})$, 
$$
v_t=D\Sigma_t(\phi)\hat{\phi}\in C,
$$
In particular, $v_0=\hat{\phi}$. The map $v$ is continuous, the restriction of $v:(-\infty,t_{\phi})\to\R^n$ to the interval $[0,t_{\phi})$ is differentiable, and
$$
v'(t)=Df(x^{\phi}_t)v_t\quad\text{for every}\quad t\in[0,t_{\phi}),
$$
with the right derivative at $t=0$.
\end{prop}

\begin{proof}
1.  Let $\phi\in U$, $\hat{\phi}\in C$, $0\le t<t_{\phi}$. For $s\le0$ with $0\le t+s$ Proposition 5.1 yields 
$$
v_t(s)=v(t+s)=(D\Sigma_{t+s}(\phi)\hat{\phi})(0)=(D\Sigma_t(\phi)\hat{\phi})(s),
$$
and for $s\le0$ with $t+s<0$,
$$
v_t(s)=v(t+s)=\hat{\phi}(t+s)=(D\Sigma_t(\phi)\hat{\phi})(s).
$$
Together, $v_t=D\Sigma_t(\phi)\hat{\phi}$. Notice that $D\Sigma_0(\phi)\hat{\phi}=\hat{\phi}$.
The fact that each segment $v_t=D\Sigma_t(\phi)\hat{\phi}$, $0\le t<t_{\phi}$, belongs to $C$ implies that $v$ is continuous.

2. Let $t>0$ with $\Omega_t\neq\emptyset$ be given. For $\phi\in\Omega_t$ consider the map 
$$
\eta^{\phi}:[0,t]\ni s\mapsto x^{\phi}(s)-\phi(0)\in\R^n.
$$
Observe that $\eta^{\phi}\in C_{0t,0}$ and 
$$
P_{0t}\phi+Z_t\eta^{\phi}=x^{\phi}|_{(-\infty,t]},
$$
which yields
$$
(P_{0t}\phi+Z_t\eta^{\phi})_s=x^{\phi}_s\in U\quad\text{for}\quad 0\le s\le t.
$$
It follows that $P_{0t}\phi+Z_t\eta^{\phi}\in dom_t$. Then
$(\eta^{\phi},\phi)$ belongs to the domain $\mathcal{O}_t$ of the map $B_t$. The map $Y_t:\Omega_t\ni\phi\mapsto \eta^{\phi}\in C_{0t,0}$ satisfies 
\begin{eqnarray*}
Y_t(\phi)(s) & = & \eta^{\phi}(s)=x^{\phi}(s)-\phi(0)\\
& = & \int_0^sf(x^{\phi}_u)du=\int_0^sf((P_{0t}\phi+Z_t\eta^{\phi})_u)du\\
& = & \int_0^sf(E_t(P_{0t}\phi+Z_tY_t(\phi),u))du=I_t(F_t(P_{0t}\phi+Z_tY_t(\phi)))(s)
\end{eqnarray*}
for all $\phi\in\Omega_t$ and $s\in[0,t]$, hence 
\begin{equation}
Y_t(\phi)=I_t(F_t(J_t(Y_t(\phi),\phi)))\qquad(=B_t(Y_t(\phi),\phi))\quad\text{for all}\quad\phi\in\Omega_t.
\end{equation}
3. Proof that the map $Y$ is $C^1_{\ast}$-smooth with
$$
v^{\phi,\hat{\phi}}(s)=(DY_t(\phi)\hat{\phi})(s)+(P_{0t}\hat{\phi})(s)\quad\text{for all}\quad s\in[0,t],\phi\in\Omega_t,\hat{\phi}\in C.
$$
By Part 2, $(Y_t(\phi),\phi)\in\mathcal{O}_t$ for all $\phi\in\Omega_t$.
With the shift map
$$
\Delta_t:C\to C_t,\quad(\Delta_t\phi)(s)=\phi(s-t),
$$
and the restriction map
$$
R_t:C_t\to C_{0t},\quad R_t\chi=\chi|_{[0,t]},
$$
which are both linear and continuous, 
$$
Y_t(\phi)=R_t(\Delta_t\circ\Sigma_t(\phi)-P_{0t}\phi)\quad\text{for all}\quad\phi\in\Omega_t.
$$
This shows that the map $Y_t$  is $C^1_{\ast}$-smooth, and for all $\phi\in\Omega_t$, $\hat{\phi}\in C$, $s\in[0,t]$,
\begin{eqnarray*}
(DY_t(\phi)\hat{\phi})(s) & = & (R_t\Delta_tD\Sigma_t(\phi)\hat{\phi})(s)-(R_tP_{0t}\hat{\phi})(s)\\
& = & (D\Sigma_t(\phi)\hat{\phi})(s-t)-\hat{\phi}(0)\\
& = & (D\Sigma_s(\phi)\hat{\phi})(0)-\hat{\phi}(0)\quad\text{(see Proposition 5.1)}\\
& = & v^{\phi,\hat{\phi}}(s)-\hat{\phi}(0)\\
& = & v^{\phi,\hat{\phi}}(s)-P_{0t}\hat{\phi}(s).
\end{eqnarray*}
For all $s\le t$ and $\phi\in\Omega_t$, $\hat{\phi}\in C$  we infer
\begin{equation}
(P_{0t}\hat{\phi})(s)+(Z_tDY_t(\phi)\hat{\phi})(s)=v^{\phi,\hat{\phi}}(s).
\end{equation}
4. Differentiation of Eq. (5.1) yields
\begin{equation}
DY_t(\phi)\hat{\phi}=I_tDF_t(J_t(Y_t(\phi),\phi))J_t(DY_t(\phi)\hat{\phi},\hat{\phi})\qquad\text{for all}\quad\phi\in\Omega_t, \hat{\phi}\in C.
\end{equation}
For such $\phi$ and $\hat{\phi}$ and for each $s\in[0,t]$,
\begin{eqnarray*}
 v^{\phi,\hat{\phi}}(s) & = & (DY_t(\phi)\hat{\phi})(s)+\hat{\phi}(0)\quad\text{(see Part 3)}\\
& = & \int_0^sDf((P_{0t}\phi)_u+(Z_tY_t(\phi))_u)((P_{0t}\hat{\phi})_u+(Z_tDY_t(\phi)\hat{\phi})_u)du\\
& & +\hat{\phi}(0)\qquad\text{(with Eq. (5.3) and  Proposition 2.2)}\\
& = & \int_0^sDf(x^{\phi}_u)v^{\phi,\hat{\phi}}_udu+\hat{\phi}(0)\quad\text{(with Eq. (5.2))}.
\end{eqnarray*}
Differentiation at $t>0$ yields 
$$
(v^{\phi,\hat{\phi}})'(t)=Df(x^{\phi}_t)v^{\phi,\hat{\phi}}_t.
$$
At $s=0$ we obtain 
$$
(v^{\phi,\hat{\phi}})'(0)=Df(x^{\phi}_0)v^{\phi,\hat{\phi}}_0=Df(\phi)\hat{\phi}
$$
with the right derivative. 
\end{proof}

\bigskip

\begin{center}
PART II
\end{center}

\bigskip

\section{Processes for nonautonomous delay differential equations}

In this section it is convenient to use the notation $C_n=C((-\infty,0],\R^n)$.
Let a set $V\subset \R\times C_n$ and a map $g:V\to\R^n$ be given. A solution of Eq. (1.3),
$$
x'(t)=g(t,x_t)
$$
on an interval $I\subset\R$ is a map $x:(-\infty,0]+I\to\R^n$ such that $(t,x_t)\in V$ for all $t\in I$ and 
the restriction $x|_I$ is differentiable and Eq. (1.3) holds for all $t\in I$ (in case $I$ has a minimum $t_0$, with the right derivative at $t_0$). For $(t_0,\phi)\in V$ a solution of the initial value problem 
\begin{equation}
x'(t)=g(t,x_t)\quad\text{for}\quad t\ge t_0,\quad x_{t_0}=\phi,
\end{equation}
is a solution $x$ of Eq. (1.3) on some interval $[t_0,t_e)$, $t_0<t_e\le\infty$, which satisfies $x_{t_0}=\phi$.

Let $p_n:C_{n+1}\to C_n$ denote the continuous linear map forgetting the first component.
For $V$ and $g$ as above define the domain
$$
U_g=\{\psi\in C_{n+1}:(\psi_1(0),p_n\psi)\in V\}
$$
and the map $f_g:C_{n+1}\supset U_g\to\R^{n+1}$ by
$$
f_g(\psi)=(1,g(\psi_1(0),p_n\psi)),
$$
so that the autonomous differential equation
\begin{equation}
y'(t)=f_g(y_t)
\end{equation}
written in components $y=(y_1,z)=(r,z)$ becomes
\begin{eqnarray*}
r'(t) & = & 1,\\
z'(t) & = & g(r(t),z_t).
\end{eqnarray*}
For $t\in\R$ given define $t_{\ast}\in C_1$ by $t_{\ast}(u)=t+u$.

\begin{prop}
(i) If $x:(-\infty,t_x)\to\R^n$ is a solution on $[t_0,t_x)$ of the IVP (6.1) then the map $(-\infty,t_x-t_0)\ni s\mapsto(s+t_0,x(s+t_0))\in \R^{n+1}$ is a solution on $[0,t_x-t_0)$ of the IVP
\begin{eqnarray}
r'(s) & = & 1\quad\text{for}\quad s\ge0,\quad r_0=t_{0\ast},\\
z'(s) & = & g(r(s),z_s)\quad\text{for}\quad s\ge0,\quad z_0=\phi.
\end{eqnarray}
(ii) If $y=(r,z)$ is a solution on $[0,t_y)$ of the IVP (6.3)-(6.4) then $x:(-\infty,t_0+t_y)\to\R^n$ given by 
$x(\tau)=z(\tau-t_0)$ is a solution on $[t_0,t_0+t_y)$ of the IVP (6.1).
\end{prop}

\begin{proof}
1. Proof of (i). Let a solution  $x:(-\infty,t_x)\to\R^n$ on $[t_0,t_x)$ of the IVP (6.1) be given and define $y:(-\infty,t_x-t_0)\to\R^{n+1}$, $y=(r,z)$ with $z:(-\infty,t_x-t_0)\to\R^n$ and $r=y_1$, by
$$
r(t)=t+t_0\quad\text{for}\quad t<t_x-t_0
$$
and
$$
z(t)=x(t+t_0)\quad\text{for all}\quad t<t_x-t_0.
$$
For $0\le t<t_x-t_0$ we get $(y_1(t),z_t)=(r(t),z_t)=(t+t_0,x_{t+t_0})\in V$. This yields $y_t\in U_g$ for $0\le t<t_x-t_0$. Obviously, $r'(s)=1$ for $0\le s<t_x-t_0$ (with the right derivative at $s=0$) and $r_0=t_{0\ast}$.
Also, for $u\le0$,
$$
z_0(u)=z(u)=x(u+t_0)=x_{t_0}(u)=\phi(u),
$$
hence $z_0=\phi$. For $0\le s<t_x-t_0$ we get 
$$
z'(s)=x'(s+t_0)=g(s+t_0,x_{s+t_0})=g(r(s),z_s),
$$
with the right derivative at $s=0$.

2. Proof of (ii). Let a solution $y=(r,z)$ on $[0,t_y)$ of the IVP (6.3)-(6.4) be given and define $x:(-\infty,t_0+t_y)\to\R^n$ by $x(\tau)=z(\tau-t_0)$.  We have $r(s)=s+t_0$ for $0\le s<t_y$. For $t_0\le\tau<t_0+t_y$ we obtain from $y_{\tau-t_0}\in U_g$ that $(\tau,x_{\tau})=(r(\tau-t_0),z_{\tau-t_0})$ belongs to $V$. Also,
$$
x'(\tau)=z'(\tau-t_0)=g(r(\tau-t_0),z_{\tau-t_0})=g(\tau,x_{\tau})
$$
(with right derivatives at $t_0$ and at $0$, respectively) while for $u\le0$,
$$
x_{t_0}(u)=x(t_0+u)=z(t_0+u-t_0)=z(u)=\phi(u).
$$
\end{proof}

Suppose now that $V$ is open and $g$ is $C^1_{\ast}$-smooth. Then $U_g\subset C_{n+1}$ is open as the preimage of $V$ under a continuous linear map, and the map $f_g:C_{n+1}\supset U_g\to\R^{n+1}$ is $C^1_{\ast}$-smooth. It follows that the solutions of the Eq. (6.2) define a continuous semiflow $\Sigma_g:[0,\infty)\times C_{n+1}\supset\Omega_g\to C_{n+1}$ on $U_g$, with all solution operators $C^1_{\ast}$-smooth. The set
$$
dom:\{(t,t_0,\phi)\in\R^2\times C_n:t_0\le t,\,\,(t-t_0,t_{0\ast},\phi)\in\Omega_g\}
$$
is an open subset of the set $\{(t,t_0)\in\R^2:t_0\le t\}\times C_n$ as it is the preimage of $\Omega_g$ under a continuous map into $[0,\infty)\times C_{n+1}$,
and the {\it process} $P:\{(t,t_0)\in\R^2:t_0\le t\}\times C_n\supset dom\to C_n$ given by
$$
P(t,t_0,\phi)=p_n\Sigma_g(t-t_0,t_{0\ast},\phi)
$$
is continuous. For every $t\ge t_0$ with $\emptyset\neq\Omega_{g,t-t_0}\subset U_g\subset C_{n+1}$ the non-empty set
$$
dom_{t,t_0}=\{\phi\in C_n:(t,t_0,\phi)\in dom\}=\{\phi\in C_n:(t_{0\ast},\phi)\in\Omega_{g,t-t_0}\}
$$
is open, and the map
$$
P(t,t_0,\cdot):C_n\supset dom_{t,t_0}\to C_n
$$
is $C^1_{\ast}$-smooth.

\begin{cor} 
(Maximal solutions, uniqueness)
For every $(t,t_0,\phi)\in\,dom$ there exists a solution $x=x^{t_0,\phi}$ of the IVP (6.1) so that any other solution of the
same IVP is a restriction of $x$, and $P(t,t_0,\phi)=x_t$.
\end{cor}

\begin{proof}
1. Let $(t,t_0,\phi)\in dom$ be given. The first part of the assertion follows from results on the autonomous Eq. (6.2), by means of Proposition 6.1.

2. We have $t_0\le t$ and $(t-t_0,t_{0\ast},\phi)\in\Omega_g$. Let $y:(-\infty,t_y)\to\R^{n+1}$ denote the maximal solution of the IVP
\begin{equation}
y'(s)=f_g(y_s)\quad\text{for}\quad s\ge0,\quad y_0=(t_{0\ast},\phi)
\end{equation}
and write $y=(r,z)$ with $r=y_1$.  Then
$$
P(t,t_0,\phi)=p_n\Sigma_g(t-t_0,t_{0\ast},\phi)=z_{t-t_0}.
$$
Proposition 6.1 says that $\tilde{x}:(-\infty, t_0+t_y)\to\R^n$ given by
$$
\tilde{x}(\tau)=z(\tau-t_0).
$$
is a solution on $[t_0,t_0+t_y)$ of the IVP (6.1). According to the first part of the assertion $\tilde{x}$ is a restriction of the maximal solution $x$ of the IVP (6.1), hence
$$
P(t,t_0,\phi)=z_{t-t_0}=\tilde{x}_t=x_t.
$$
\end{proof}
 
\begin{cor}
For all $(t_0,\phi)\in V$, $(t_0,t_0,\phi)\in dom$ and $P(t_0,t_0,\phi)=\phi$, and for all $t_0\le t\le s$ with $(t,t_0,\phi)\in dom$ and $(s,t,P(t,t_0,\phi))\in dom$,
$$
(s,t_0,\phi)\in dom\quad\text{and}\quad P(s,t_0,\phi)=P(s,t,P(t,t_0,\phi)).
$$
\end{cor}

\begin{proof}
1. For $(t_0,\phi)\in V$ we have $(t_{0\ast},\phi)\in U_g$, hence $(0,t_{0\ast},\phi)\in\Omega_g$. It follows that
$(t_0,t_0,\phi)\in dom$ and $P(t_0,t_0,\phi)=p_n\Sigma_g(0,t_{0\ast},\phi)=\phi$.

2. Suppose  $t_0\le t\le s$, $(t,t_0,\phi)\in dom$, and $(s,t,P(t,t_0,\phi))\in dom$. Then $(t-t_0,t_{0\ast},\phi)\in\Omega_g$. As the solution of the IVP (6.3) is given by $r(s)=s+t_0$ we see that the first component $r_{t-t_0}$ of $\Sigma_g(t-t_0,t_{0\ast},\phi)$ satisfies 
$$
r_{t-t_0}(u)=r(t-t_0+u)=(t-t_0+u)+t_0=t+u=t_{\ast}(u)\quad\text{for all}\quad u\le0,
$$
or, $r_{t-t_0}=t_{\ast}$. It follows that 
$$
\Sigma_g(t-t_0,t_{0\ast},\phi)=(t_{\ast},P(t,t_0,\phi)).
$$
Using this and $(s,t,P(t,t_0,\phi))\in dom$ we get
$$
(s-t,\Sigma_g(t-t_0,t_{0\ast},\phi))=(s-t,t_{\ast},P(t,t_0,\phi))\in\Omega_g.
$$
Now the properties of $\Sigma_g$ yield
$$
(s-t_0,t_{0\ast},\phi)=((s-t)+(t-t_0),t_{0\ast},\phi)\in\Omega_g,
$$
hence $(s,t_0,\phi)\in dom$ and 
$$
P(s,t_0,\phi)=p_n\Sigma_g(s-t_0,t_{0\ast},\phi)=p_n\Sigma(s-t,\Sigma(t-t_0,t_{0\ast},\phi))
$$
$$
=p_n\Sigma_g(s-t,t_{\ast},P(t,t_0,\phi))=P(s,t,P(t,t_0,\phi)).
$$
\end{proof}

\section{Volterra integro-differential equations}

Consider the Volterra integro-differential equation (1.4),
$$
x'(t)=\int_0^tk(t,s)h(x(s))ds
$$
with $k:\R^2\to\R^{n\times n}$ and $h:\R^n\to\R^n$ continuously differentiable. With $K:\R^2\to\R^{n\times n}$ given by $K(t,s)=k(t,t+s)$ we may write 
\begin{eqnarray}
x'(t) & = & \int^0_{-t}k(t,t+s)h(x(t+s))ds\nonumber\\
& = & \int_{-t}^0K(t,s)h(x_t(s))ds,
\end{eqnarray}
with $x_t\in C$. A solution of Eq. (7.1) would be a continuous map $x:(-\infty,t_e)\to\R^n$, $0<t_e\le\infty$, whose restriction to the interval $(0,t_e)$ is differentiable and satisfies Eq. (7.1). We look for a map $g:\R\times C\to\R^n$ so that every solution of Eq. (7.1) also is a solution on $(0,t_e)$ of  Eq. (1.3) ,
$$
x'(t)=g(t,x_t).
$$
In order to avoid advanced arguments of segments of solutions we employ the {\it odd prolongation map}
$$
P_o:C\to C_{\infty}
$$
given by $P_o\phi(s)=\phi(s)$ for $s\le0$ and $P\phi(s)=2\phi(0)-\phi(-s)$ for $0<s$. The map $P_o$ is linear and continuous. (We could also use constant prolongation for the present purpose. Odd prolongation has the advantage that it defines a continuous linear map $C^1\to C^1_{\infty}$. This plays a role when considering nonautonomous equations with discrete delay, like the {\it pantograph equation}
$$
x'(t)=a\,x(\lambda t)+b\,x(t)
$$
with $0<\lambda<1$.)

Next, consider the  substitution operator 
$$ 
S_H:C_{\infty}\ni\phi\mapsto H\circ\phi\in C_{\infty}
$$
which is defined for every continuous map $H:\R^n\to\R^n$,
and the linear integration operator
$$
I:C_{\infty}\to C^1_{\infty}
$$
given by  $(I\psi)(u)=\int_{-u}^0K(u,s)\psi(s)ds$, and the operator 
$$
J:C^1_{\infty}\times\R\to\R^n
$$
given by
$$
J(\psi,t)=\int_{-t}^0K(t,s)\psi(s)ds=Ev_{\infty,1}(I\psi,t).
$$ 
Define $g:\R\times C\to\R^n$ by
$$
g(t,\phi)=J((S_h\circ P_o)(\phi),t)=\int_{-t}^0K(t,s)h((P_o\phi)(s))ds
$$
and observe that indeed for every solution $x:(-\infty,t_e)\to\R^n$ of Eq. (7.1) and for all $t\in(0,t_e)$ we have
$$
g(t,x_t)=\int_{-t}^0K(t,s)h((P_ox_t)(s))ds=\int_{-t}^0K(t,s)h(x_t(s))ds.
$$
In order to show that the map $g$ is $C^1_F$-smooth recall first that the evaluation map $Ev_{\infty,1}$ is $C^1_F$-smooth. Therefore the map $J$ is $C^1_F$-smooth provided the linear map $I$ is continuous. It follows that the map $g$ is $C^1_F$-smooth provided $I$ is continuous and $S_h$ is $C^1_F$-smooth. The next propositions establish these remaining smoothness properties.

\begin{prop}
The linear map $I$ is continuous.
\end{prop}  

\begin{proof}
Use the relations
\begin{eqnarray*}
|I\psi|_{\infty,j} & = & \max_{-j\le u\le j}\left|\int_{-u}^0K(u,s)\psi(s)ds\right|\le|j|\max_{-j\le u\le j,-j\le s\le j}|K(u,s)||\psi|_{\infty,j},\\
(I\psi)'(u) & = & -K(u,u)\psi(u)+\int_{-u}^0\partial_1K(u,s)\psi(s)ds,\\
|(I\psi)'|_{\infty,j} & \le & \max_{-j\le u\le j}|K(u,u)||\psi|_{\infty,j}+|j|\max_{-j\le u\le j,-j\le s\le j}|\partial_1K(u,s)||\psi|_{\infty,j}
\end{eqnarray*}
for all $j\in\N$ and $\psi\in C_{\infty}$.
\end{proof}

\begin{prop}
If $H:\R^n\to\R^n$ is continuous then the map $S_H$ is continuous. In case $H$ is continuously differentiable the map $S_H$ is $C^1_F$-smooth, with
$$
(DS_H(\phi)\chi)(t)=DH(\phi(t))\chi(t).
$$
\end{prop}

\begin{proof}
1. For $j\in\N$ set
$$
N_j=\left\{\phi\in C_{\infty}:|\phi|_{\infty,j}<\frac{1}{j}\right\}.
$$
Let $H$ be continuous. Let $\phi\in C_{\infty}$. For continuity of $S_H$ at $\phi$ we need that for every $j\in\N$ there exists $k\in \N$ such that for all
$\chi\in C_{\infty}$ with $\chi\in\phi+N_k$ we have $S_H(\chi)\in S_ H(\phi)+N_j$.  Let $j\in\N$ be given. Choose a compact neighbourhood $W$ of $\phi([-j,j])$. As $H$ is uniformly continuous on $W$ there exists $\delta>0$ with $|H(y)-H(x)|<\frac{1}{j}$ for all $x,y$ in $W$ with $|y-x|<\delta$. Choose $k\in\N$ with $k\ge j$ and $\frac{1}{k}<\delta$ and $\chi([-j,j])\subset W$ for all $\chi\in C_{\infty}$ with $|\chi-\phi|_{\infty,k}<\frac{1}{k}$ (or equivalently, $\chi\in\phi+N_k$). For such $\chi$ and for all $s\in[-j,j]$ we get 
$$
|H(\chi(s))-H(\phi(s))|<\frac{1}{j},
$$ 
hence $S_H(\chi)\in S_H(\phi)+N_j$.

2.  Let $H$ be continuously differentiable. 

2.1. (Existence of directional derivatives) Let $\phi\in C_{\infty}$ and $\chi\in C_{\infty}$be given. Define $A(\phi,\chi)\in C_{\infty}$ by $A(\phi,\chi)(s)=DH(\phi(s))\chi(s)$. It is sufficient to show that for every $j\in\N$ we have
$$
|t^{-1}(S_H(\phi+t\chi)-S_H(\phi))-A(\phi,\chi)|_{\infty,j}\to0\quad\text{as}\quad0\neq t\to0.
$$ 
Let $j\in\N$ be given. For all
reals $t\neq0$,
$$
|t^{-1}(S_H(\phi+t\chi)-S_H(\phi))-A(\phi,\chi)|_{\infty,j}
$$
$$
=\max_{-j\le s\le j}|t^{-1}(H(\phi(s)+t\chi(s))-H(\phi(s)))-DH(\phi(s))\chi(s)|
$$
$$
=\max_{-j\le s\le j}\left|\int_0^1(DH(\phi(s)+ut\chi(s))\chi(s)-DH(\phi(s))\chi(s))du\right|
$$
$$
\le\max_{-j\le s\le j}\max_{|v|\le|t|}|DH(\phi(s)+v\chi(s))-DH(\phi(s)||\chi(s)|
$$
$$
\le|\chi|_{\infty,j}\max_{-j\le s\le j}\max_{|v|\le|t|}|DH(\phi(s)+v\chi(s))-DH(\phi(s)|.
$$
As $DH$ is continuous, arguments as in Part 1 of the proof  can be used in order to deduce from the previous estimate that we have
$$
\lim_{0\neq t\to0}|t^{-1}(S_H(\phi+t\chi)-S_H(\phi))-A(\phi,\chi)|_{\infty,j}=0.
$$

2.2. Each map $DS_H(\phi):C_{\infty}\ni\chi\mapsto A(\phi,\chi)\in C_{\infty}$, $\phi\in C_{\infty}$, is linear. Continuity follows from the estimates $|A(\phi,\chi)|_{\infty,j}\le\max_{-j\le s\le j}|DH(\phi(s))||\chi|_{\infty,j}$ for all $j\in\N$ and all $\chi\in C_{\infty}$.

2.3. It remains to show that $DS_H:C_{\infty}\ni\phi\mapsto DS_H(\phi)\in L_c(C_{\infty},C_{\infty})$ is continuous with respect to the topology $\beta$ on $L_c(C_{\infty},C_{\infty})$. Let $B\subset C_{\infty}$ be bounded and let $j\in\N$. According to Remark 8.1 
(\cite[Remark 2.1 (iii)]{W9}) we have to find an integer $k\ge j$ such that
$$
DS_H(\psi)\chi-DS_H(\phi)\chi=A(\psi,\chi)-A(\phi,\chi)\in N_j\quad\text{for all}\quad\psi\in\phi+ N_k\quad\text{and}\quad\chi\in B.
$$
By \cite[Theorem 1.37]{R}, $b_j=\sup_{\chi\in B}|\chi|_{\infty,j}<\infty$. Using arguments as in Part 1 of the proof one finds an integer $k\ge j$ such that for all $\psi\in\phi+N_k$ we have
$$
\max_{-j\le s\le j}|DH(\psi(s))-DH(\phi(s))|<\frac{1}{jb_j}.
$$
For such $\psi$ and for all $\chi\in B$ we infer
\begin{eqnarray*}
|A(\psi,\chi)-A(\phi,\chi)|_{\infty,j} & \le & \max_{-j\le s\le j}|DH(\psi(s))-DH(\phi(s))||\chi|_{\infty,j}\\
& \le & \max_{-j\le s\le j}|DH(\psi(s))-DH(\phi(s))|b_j<\frac{1}{j},
\end{eqnarray*}
or, 
$$
A(\psi,\chi)-A(\phi,\chi)\in N_j\quad\text{for all}\quad\psi\in\phi+ N_k\quad\text{and}\quad\chi\in B.
$$
\end{proof}

\begin{cor}
The map $g:\R\times C\to\R^n$ given by $g(t,\phi)=J((S_h\circ P_o)(\phi),t)$ is $C^1_F$-smooth.
\end{cor}

The results of the previous section apply to the nonautonomous equation (1.3) with $g$ from the preceding corollary and yield a continuous process of solution operators $P(t,t_0)$ which are defined on open subsets of $C$, and which are $C^1_F$-smooth. 

\section{Appendix: Uniform convergence of continuous linear maps on bounded subsets, $C^1_F$-smoothness}

Let $V,W$ be topological vector spaces over $\R$ or $\C$. On $L_c=L_c(V,W)$ the topology $\beta$  of uniform convergence on bounded sets is defined as follows.
For a neighbourhood $N$ of $0$ in $W$ and a bounded set $B\subset V$ the neighbourhood $U_{N,B}$ of $0$ in $L_c$ is defined as
$$
U_{N,B}=\{A\in L_c:TB\subset N\}.
$$
Every finite intersection of such sets $U_{N_j,B_j}$, $j\in\{1,\ldots,J\}$, contains a set of the same kind, because of the inclusion
$$
\cap_{j=1}^JU_{N_j,B_j}\supset\{T\in L_c:T(\cup_{j=1}^JB_j)\subset\cap_{j=1}^JN_j\}
$$
and since finite unions of bounded sets are bounded and finite intersections of neighbourhoods of $0$ are neighbourhoods of $0$. Then the topology $\beta$ is the set of all $O\subset L_c$ which have the property that for each $A\in O$ there exist a neighbourhood $N$ of $0$ in $W$ and a bounded set $B\subset V$ with
$A+U_{N,B}\subset O$.

\medskip

We call a map $A$ from a topological space $T$ into $L_c$ $\beta$-continuous at a point $t\in T$ if it is continuous at $t$ with respect to the topology $\beta$ on $L_c$.

\begin{remark} \cite[Remark 2.1 (iii)]{W9} In order to verify $\beta$-continuity of a map $A:T\to L_c$, $T$ a topological space, at some $t\in T$ one has to show that, given a bounded subset $B\subset V$ and a neighbourhood $N$ of $0$ in $W$, there exists a neighbourhood $N_t$ of $t$ in $T$  such that for all $s\in N_t$ we have $(A(s)-A(t))(B)\subset N$. 

\medskip

In case $T$ has countable neighbourhood bases 
the map $A$ is $\beta$-continuous  at $t\in T$ if and only if for any sequence $T\ni t_j \to t$ we have $A(t_j)\to A(t)$. For $A(t_j)\to A(t)$ we need that given a bounded subset $B\subset V$ and a neighbourhood $N$ of $0$ in $W$, there exists $J\in\N$ with
$$
(A(t_j)-A(t))(B)\subset N\quad\text{for all integers}\quad j\ge J.
$$
\end{remark}

\begin{prop} \cite[Proposition 3.2]{W9}
Let $F$ and $G$ be Fr\'echet spaces, $U\subset F$ open. A map $g:U\to G$ is $C^1_F$-smooth if and only if it is  $C^1_{MB}$-smooth with $U\ni u\mapsto Dg(u)\in L_c(F,G)$ $\beta$-continuous.
\end{prop}

Continuous linear maps $L:F\to G$ between Fr\'echet spaces are $C^1_F$-smooth since they are $C^1_{MB}$-smooth with constant derivative $DL(u)=L$ for all $u\in F$, and differentiation $g\mapsto Dg$ of $C^1_F$-maps $U\to G$ is linear.

The following two propositions are included for convenience, without being used in Sections 2-7.

\begin{prop} \cite[Proposition 3.3]{W9}
In case $E$ is a finite-dimensional normed space each $C^1_{MB}$-map $g:E\supset U\to G$ is $C^1_F$-smooth.
\end{prop}

\cite[Section 8]{W9} contains examples of maps on infinite-dimensional Banach spaces which are $C^1_{MB}$-smooth but not $C^1_F$-smooth.

\begin{prop} \cite[Proposition 3.4]{W9}
For Banach spaces $F$ and $G$ and $U\subset F$ open a map $g:F\supset U\to G$ is $C^1_F$-smooth if and only if there exists a continuous map $D_g:U\to L_c(F,G)$  such that for every $u\in U$ and\\
(F) $\quad$ for every $\epsilon>0$ there exists $\delta>0$ with
$$
|g(v)-g(u)-D_g(u)(v-u)|\le\epsilon|v-u|\quad\text{for all}\quad v\in U\quad\text{with}\quad |v-u|<\delta.
$$ 
In this case, $D_g(u)v$ is the directional derivative $Dg(u)v$, for every $u\in U,v\in F$. 
\end{prop}

\begin{prop} \cite[Proposition 3.5]{W9}
(Chain rule). If $g:F\supset U\to G$ and $h:G\supset V\to H$ are $C^1_F$-maps, with $g(U)\subset V$, then also
$h\circ g$ is a $C^1_F$-map.
\end{prop}

\begin{prop} \cite[Proposition 3.6]{W9} 
Let Fr\'echet spaces $F_1,F_2,G$ be given. For a continuous map $g:F_1\times F_2\supset U\to G$, $U$ open, the following statements are equivalent.\\
(i) For all $(u_1,u_2)\in U$ and all $v_k\in F_k$, $k\in\{1,2\}$, $g$ has a partial derivative $D_kg(u_1,u_2)v_k\in G$, all maps
$$
D_kg(u_1,u_2):F_k\to G,\quad(u_1,u_2)\in U,\quad k\in\{1,2\},
$$
are linear and continuous, and the maps
$$
U\ni(u_1,u_2)\mapsto D_kg(u_1,u_2)\in L_c(F_k,G),\quad k\in\{1,2\},
$$ 
are $\beta$-continuous.\\
(ii) $g$ is $C^1_F$-smooth.\\
In this case,
$$
Dg(u_1,u_2)(v_1,v_2)=D_1g(u_1,u_2)v_1+D_2g(u_1,u_2)v_2
$$
for all $(u_1,u_2)\in U$, $v_1\in F_1$, $v_2\in F_2$.
\end{prop}

\begin{thm} \cite[Theorem 5.2]{W9}
Let a Fr\'echet space $T$, a Banach space $B$, open sets $V\subset T$ and $O_B\subset B$, and a $C^1_F$-map $A:V\times O_B\to B$ be given. Assume that for a closed set $M\subset O_B$ we have $A(V\times M)\subset M$, and $A$ is a uniform contraction  in the sense that there exists $k\in[0,1)$ so that 
$$
|A(t,x)-A(t,y)|\le k|x-y|
$$
for all $t\in V,x\in O_B,y\in O_B$. Then the map $g:V\to B$ given by $g(t)=A(t,g(t))\in M$ is $C^1_F$-smooth.
\end{thm}

\begin{prop} (See \cite[Proposition 9.1 (iii)]{W9} for $T=\infty$)
The map 
$$
E^{10}_{\infty}:C^1_{\infty}\times\R\to C,\quad E^{10}_{\infty}(\phi,t)=\phi_t,
$$ 
is $C^1_F$-smooth, with
\begin{eqnarray*}
D_1E^{10}_{\infty}(\phi,t)\hat{\phi} & = & \hat{\phi}_t\quad\text{and}\\
D_2E^{10}_{\infty}(\phi,t)t_{\ast} & = & t_{\ast}(\phi')_t.
\end{eqnarray*}
\end{prop}

\end{document}